\documentclass[12pt]{article}
\font\smallit=cmti10

\usepackage{amssymb,latexsym}

\usepackage{amsmath}
\usepackage[colorlinks=true, pdfstartview=FitV, linkcolor=blue, citecolor=blue, urlcolor=red]{hyperref} 
\usepackage[arrow,matrix,curve,frame]{xy}
\usepackage{graphicx}

\makeatletter

\renewcommand\section{\@startsection {section}{1}{\z@}
{-30pt \@plus -1ex \@minus -.2ex}
{2.3ex \@plus.2ex}
{\normalfont\normalsize\bfseries}}

\renewcommand\subsection{\@startsection{subsection}{2}{\z@}
{-3.25ex\@plus -1ex \@minus -.2ex}
{1.5ex \@plus .2ex}
{\normalfont\normalsize\bfseries}}

\renewcommand{\@seccntformat}[1]{\csname the#1\endcsname. }

\makeatother

\newtheorem{thm}{Theorem}[section]

\newcommand{\ABC}{\ensuremath{\mathcal{A}}}
\newcommand{\NN}{\ensuremath{\mathbb{N}}}

\begin{document}

\begin{center}
\uppercase{\bf More Kolakoski Sequences}
\vskip 20pt
{\bf Bernd Sing}\\
{\smallit Department~of~Computer~Science,~Mathematics~\&~Physics, University~of~the~West~Indies, Cave~Hill, P.O.~Box~64, Bridgetown,~BB11000, Barbados, West~Indies}\\
{\tt bernd.sing@cavehill.uwi.edu}
\end{center}
\vskip 30pt

\centerline{\bf Abstract}
\noindent
Our goal in this article is to review the known properties of the mysterious Kolakoski sequence and at the same time look at generalizations of it over arbitrary two letter alphabets. Our primary focus will here be the case where one of the letters is odd while the other is even, since in the other cases the sequences in question can be rewritten as (well-known) primitive substitution sequences. We will look at word and letter frequencies, squares, palindromes and complexity.


\thispagestyle{empty} 
\baselineskip=15pt 
\vskip 30pt

\section{Introduction}

A one-sided infinite sequence $z$ over the alphabet $\ABC=\{1,2\}$ is
called a (classical) \emph{Kolakoski sequence}, if it equals the sequence defined by
its run-lengths, i.e.:
\begin{equation*}
\begin{array}{ccccccccccccccc}
z & = & \underbrace{22} & \underbrace{11} & \underbrace{2} &
\underbrace{1} & \underbrace{22} & \underbrace{1} & \underbrace{22} &
\underbrace{11} & \underbrace{2} & \underbrace{11} & \ldots && \\
&& 2 & 2 & 1 & 1 & 2 & 1 & 2 & 2 & 1 & 2 & \ldots & = & z.
\end{array}
\end{equation*}
Here, a \emph{run} is a maximal subword consisting of identical letters. The
sequence $z'=1z$ is the only other sequence which has this property.

This sequence was introduced by Kolakoski in~\cite{Kol65} who asked ``What is the $n$th term? Is the sequence periodic?''\footnote{
	The first question is still studied today, see \cite{Ste06} and \cite{FF10}. In these articles, recursive formulae for the $n$th term are derived thus answering the first question.
}
This sequence has attracted attention over the years since, although it is easy to define, it resists any attempt to reveal even some of its most basic properties like recurrence or the frequency of its letters. There is even some prize money offered for answering some of these question about its properties, see \cite{Kim79,Kim3W}. The maybe most basic question is known as \emph{Keane's question} \cite{Kea91}:
\begin{quote}
Does the frequency of the symbol $1$ in $z=221121\ldots$ exist, and is it equal to $\frac12$?
\end{quote}
The line of attack in trying to prove this question has often been to detect some structure by rewriting the generation rule of the Kolakoski sequence in some sort of generalized substitution rule, see for example \cite{Ste96,Dek97},  \cite[Section 4.4]{Fogg} and references therein. However, these attempts have not been successful in answering Keane's question.

Our goal in this article is more humble -- we want to give an overview of (the little) what is known about the Kolakoski sequence but at the same time look at generalizations to arbitrary two-letter alphabets $\ABC=\{r,s\}$ where $r$ and $s$ are natural numbers (with $r\neq s$). We note that ``$10$'' or ``$439$'' in this generalization is \textit{one} letter not two or three, and we can (well, if we really want to) examine the Kolakoski sequence(s) over the alphabet $\ABC=\{10,439\}$. 

If we do this generalization, then we find some easy cases for which we can answer Keane's question immediately. For this we use the observation made in \cite{CKL92}: One can obtain the Kolakoski sequence $z$ above by starting with $2$ as a seed and iterating the two substitutions 
\begin{equation*}
\sigma^{}_0: \begin{array}{lcl} 1 & \mapsto & 1 \\ 2 & \mapsto & 11,\end{array} 
\quad \text{and} \quad
\sigma^{}_1: \begin{array}{lcl} 1 & \mapsto & 2 \\ 2 & \mapsto & 22 \end{array}
\end{equation*}
alternatingly, i.e., $\sigma^{}_0$ substitutes letters on even positions and
$\sigma^{}_1$ letters on odd positions:
\begin{equation*}
2 \mapsto 22 \mapsto 2211 \mapsto 221121 \mapsto 221121221 \mapsto \ldots
\end{equation*}
Clearly, the iterates converge to the Kolakoski sequence $z$ (in the
obvious product topology), and $z$ is the unique (one-sided) fixed point
of this iteration.  

Similarly, a (generalized) Kolakoski sequence over an alphabet $\ABC=\{r,s\}$, which is again also equal to
the sequence of its run-lengths, can be obtained by iterating the two substitutions 
\begin{equation*}
\sigma^{}_0: \begin{array}{lcl} r & \mapsto & r^r \\ s & \mapsto & r^s
\end{array} \quad \text{and} \quad
\sigma^{}_1: \begin{array}{lcl} r & \mapsto & s^r \\ s & \mapsto & s^s
\end{array} 
\end{equation*}
alternatingly. Here, $a^b$ denotes a run of $b$ $a$s, i.e., $a^b=a\ldots a$ ($b$ times).

Let us now assume that both $r$ and $s$ are even number. Building blocks of two letters $A=rr$ and $B=ss$ and applying the alternating substitution rule to them, one actually obtains a usual substitution rule for $A$ and $B$:
\begin{equation*}
\sigma: \begin{array}{lcl} A & \mapsto & A_{}^{m} B_{}^{m} \\ B & \mapsto &
  A_{}^{n} B_{}^{n} \end{array}
\end{equation*}
where $m=\frac{r}2$ and $n=\frac{s}2$. In fact, from this (primitive) substitution rule it is easy to see that the frequency of the letters $r$ and $s$ in the original sequence must be equal, see \cite{diplom,Sin03}. 

Let us now assume that both $r$ and $s$ are odd numbers. Again, building blocks of two letters helps, although we need three such blocks here: $A=rr$, $B=rs$ and $C=ss$. For these three letters one again obatins a usual (primitive) substitution rule:
\begin{equation}\label{eq:PV}
\sigma: \begin{array}{lcl} A & \mapsto & A_{}^{m} BC_{}^{m} \\ B & \mapsto &
  A_{}^{m} BC_{}^{n} \\  C & \mapsto &  A_{}^{n} BC_{}^{n} \end{array}
\end{equation}
where $m=\frac{r-1}2$ and $n=\frac{s-1}2$. From this representation it is straightforward to calculate the letter frequencies in the corresponding Kolakoski sequence. However, here the frequencies of $r$ and $s$ are not equal\footnote
{
	One can show that the substitution in \eqref{eq:PV} is a Pisot substitution with cubic Pisot-Vijaraghavan number if $2(r+s)\ge (r-s)^2$. It is a unimodular Pisot substitution if $r=s\pm 2$. In the case $2(r+s)<(r-s)^2$, all roots of the corresponding substitution matrix are greater than $1$ in modulus (and cubic algebraic numbers).  A formula for the letter frequencies in the case that one of the odd numbers is $1$ can be found in \cite{BJP08}.
}, 
see \cite{BS04,diplom}.

We will therefore look at the generalizations of the Kolakoski sequence in this article where one of the letters in the alphabet is odd while the other is even. We will not look at generalizations to three-letter alphabets (see for example   \cite{BBC05}) since there the situation is in general\footnote
{
	Of course, there are also simple cases where we can rewrite everything using one substitution rule: If the three letters are equal modulo $3$, building blocks of three letters is the key. At least, if we alternate the three letters periodically in the original sequence.	
} 
certainly worse than for two-letter (where we only alternate between two letters).

\section{Derivatives and Primitives}

Broadly speaking, there are (currently) two approaches to study a Kolakoski sequence: Either one tries to examine the set of all (infinite) sequences over $\ABC=\{r,s\}$ with the property that their run-length sequence is also a sequence over the same alphabet $\ABC=\{r,s\}$ (and the run-length sequence of the run-length sequence -- and so on -- is also a sequence over $\ABC=\{r,s\}$). Or, one tries to study the set of all possible (finite) subwords (or factors) of the Kolakoski sequence. This leads to the study of so-called $C^{\infty}$-words. We will introduce $C^{\infty}$-words in this and the next section, and will show how the former approach via sequences is used in Section~\ref{sec:Chvatal}.

We start with some basic definitions. Let $\ABC$ be an alphabet, which throughout this article will always be a two-letter alphabet $\ABC=\{r,s\}$ where $r,s\in\NN$. Then $z\in\ABC^{\NN}$ is a (one-sided infinite) \emph{sequence} of letters in $\ABC$. Any $w=w_1w_2\ldots w_n\in \ABC^n$ where $n\in\NN$ is a \emph{word} of length $n$ and we use the notation $|w|=n$ to denote the length of $w$. We denote the \emph{empty word} by $\varepsilon$. Furthermore, we use the notation $|w|_r$ and $|w|_s$ for the number of $r$s and $s$s in the word $w$, and, moreover, $|w|_v$ for the number of occurences of the word $v$ in the word $w$. 

Since we are working in a two-letter alphabet, we can define the following two properties: Let $\tilde{\cdot}$ be the operation that exchanges letters, i.e., $\tilde{r}=s$ and $\tilde{s}=r$ extended to any word $w=w_1w_2\ldots w_n$ by $\tilde{w}=\tilde{w}_1\tilde{w}_2\ldots\tilde{w}_n$. Then a sequence $z$ is called \emph{mirror invariant} if
\begin{equation*}
w \mbox{ occurs in } z \qquad \Longleftrightarrow \qquad \tilde{w} \mbox{ occurs in } z.
\end{equation*}
Similarly, the operation $\stackrel{\leftarrow}{\cdot}$ denotes the reversed word $\stackrel{\leftarrow}{w}=w_nw_{n-1}\ldots w_2w_1$ of a word $w=w_1w_2\ldots w_{n-1}w_n$, and we say that a sequence $z$ is \emph{reversal invariant} if
\begin{equation*}
w \mbox{ occurs in } z \qquad \Longleftrightarrow \qquad \stackrel{\leftarrow}{w} \mbox{ occurs in } z.
\end{equation*}

Let us have a look back at those Kolakoski sequences where either $r$ and $s$ are both odd or both even: 
\begin{itemize}
\item If $r$ and $s$ are both even, the corresponding Kolakoski sequence is not mirror invariant: If the sequence starts with $r$, then a run of $r$s is always followed by an unequally long run of $s$s, but a run of $s$s is followed by a run of $s$s of length $r$ or $s$. E.g., the sequence $z=2244222244442\ldots$ has the subword $4444224$ but not $2222442$.\\
This example also shows that such a Kolakoski sequence is not reversal invariant (e.g., $2442222$ does not appear in the previous $z$).
\item If $r$ and $s$ are both odd, the corresponding Kolakoski sequence is also not mirror invariant: We have that $r^r$ is followed by either $s^s$ or $s^r$ (and $s^s$ by either $r^s$ or $r^r$) while $r^s$ is only followed by $s^s$ (and $s^r$ by $r^r$). E.g., $313331$ appears in the Kolakoski sequence $z=3331113331313331\ldots$ while $131113$ does not.
\item However, if $r$ and $s$ are both even, the corresponding Kolakoski sequence is reversal invariant: One can extend any subword $w$ to right, say $uw$, such that $uw$ is a \emph{palindrome} (i.e., such that $uw= \stackrel{\longleftarrow}{uw}$) and a subword of the Kolakoski sequence. E.g., in the previous example, since $313331$ appears, so does $13331313331$ and this establishes immediately that $133313$ appears in the Kolakoski sequence over $\{1,3\}$. The reader may convince herself that this construction always works (one can start by the observation that $r^r$ is preceeded by either $s^s$ or $s^r$, while $r^s$ must be preceeded by $s^s$). 
\end{itemize}
 
Our goal is now to see if we can say something about these properties in the case where one of the letters $\{r,s\}$ is odd while the other is even. Here, we will closely follow \cite[Section 3]{Dek97}. From now on, we will use the following convention:
\begin{equation*}
r=\min\{r,s\} \qquad \mbox{and} \qquad s=\max\{r,s\}.
\end{equation*}
Let $w$ be a word over $\ABC=\{r,s\}$. We define the following ``\emph{differentiation}'' rule for $w$: The derivative $D(w)$ of $w$ is, in principle, the run-length sequence of $w$ except for (possibly) the first and last symbol. If $w$ is a single run of length less than $s$, we set $D(w)=\varepsilon$. If $w$ consists of more than one run and the first (last) run of $w$ is of length less than or equal to $r$, we discard this run it. If $w$ consists of more than one run and the first (last) run of $w$ has length between $r+1$ and $s$, we extend it to a run of length $s$. The word $D(w)$ is now the run-length sequence of this altered word and might be the empty word $\varepsilon$ (we use the convention $D(\varepsilon)=\varepsilon$). We say that $w$ is \emph{differentiable} if $D(w)$ is again a word over the same alphabet $\ABC=\{r,s\}$. Let us look at some examples using the alphabet $\{2,5\}$:
\begin{equation*}
\begin{array}{lllll}
D(255555222) = 55 & D(2555552) = 5 &  D(2255) = \varepsilon & D(222555)=55 & D(2222) = \varepsilon \\
D(25555552) = 6 & D(25252) = 111 & D(222522) = 51 & D(2555222) = 35 & \\
\end{array}
\end{equation*}    
Note that the words in the second line are not differentiable!

The definition of differentiable is chosen such that every subword of a Kolakoski sequence is differentiable. In fact, every subword of a Kolakoski sequence is \emph{smooth} or a \emph{$C^{\infty}$-word} with respect to this differentiation rule over the respective alphabet, i.e., it is arbitrarily often differentiable. 

We say that a word $v$ is a \emph{primitive} of a word $w$ if $D(v)=w$. From our differential rule (discarding and/or extending the first and last run), one can conclude that each (nonempty) word has at least $2r^2$ and at most $2s^2$ primitives (the factor $2$ appears since we have $D(v)=w=D(\tilde{v})$, i.e., a word and its mirrored word have the same derivative). E.g., over the alphabet $\{2,3\}$ the primitives of $33$ are: 
$\,222333$, $\,3222333$, $\,33222333$, $\,2223332$, $\,22233322$, $\,32223332$, $\,332223332$, $\,322233322$, $\,3322233322$, $\,333222$, $\,2333222$, $\,22333222$,  $\,3332223$, $\,33322233$, $\,23332223$, $\,223332223$, $\,233322233$, $\,2233322233$.

One can now use the differentiation rule to prove the following statements:

\begin{thm}
\begin{enumerate}
\item\label{it:prop:1} Kolakoski sequences are not eventually periodic (where a sequence $z$ is called \emph{eventually periodic} if there exist $m,q\in\NN$ such that $z_{i+1}\ldots z_{i+q}=z_{i+q+1}\ldots z_{i+2\,q}$ for all $i\ge m$).
\item\label{it:prop:2} For a Kolakoski sequence, mirror invariance implies recurrence (where a sequence $z$ is called recurrent if any word that occurs in $z$ does so infinitely often).
\item\label{it:prop:3} For a Kolakoski sequence, mirror invariance holds iff each $C^{\infty}$-word occurs in it.
\end{enumerate}
\end{thm}

\noindent
{\it Proof.} 
\begin{enumerate}
\item[\ref{it:prop:1}] Compare \cite{Kol65} and \cite[Example 4]{Dek79}. The reason is that a (minimal) period of length $q$ in a sequence $z$ yields a period of length $q'<q$ in its run-length sequence. Thus such a sequence $z$ cannot be equal to its run-length sequence.
\item[\ref{it:prop:2}] The proof of \cite[Proposition 3.1]{Dek97} also applies here.
\item[\ref{it:prop:3}] The proof of \cite[Proposition 2]{Dek80} also applies here.\hfill $\Box$
\end{enumerate}

We have seen above that in the case where the letters $\{r,s\}$ are both even or odd, the corresponding Kolakoski sequence is not mirror invariant. Of course, since they can be constructed using a primitive substitution rule, they are recurrent and even \emph{repetitive} (or \emph{uniformly recurrent}): Every word that occurs in the sequence does so with bounded gaps.

However, for all Kolakoski sequence over one even and one odd symbol, nothing seems to be known beyond the above implications. We don't know whether or not all $C^{\infty}$-words occur in such a Kolakoski sequence, or whether or not it is recurrent. In fact, it is even not known whether or not a Kolakoski sequence is repetitive. The problem with the last property is, of course, that the gap might be quite large, thus one has to be careful with claims based on numerical studies (as in \cite[Section 4.1.4]{Lad99}). But one can use $C^{\infty}$-words to answer the following question\footnote{
	For the question ``Given $|v|\le n$, what is the maximal possible length of $w$ such that $wvw$ is a $C^{\infty}$-word?'' see \cite[Proposition 7]{Car94}: Based on the computations in \cite{Chv94}, this length is bounded $O(n^{1.002})$, and it is conjectured to be $O(n)$. Also see \cite[Section 6.3]{DAldiss} and \cite{CDA09} on this question and its connection to Keane's question.
}: 
Given a word $w$, what is the maximal possible length of $v$ such that $wvw$ is a $C^{\infty}$-word and $w$ is not a subword of $v$? For the classical Kolakoski sequence over $\{1,2\}$ one obtains the following table:
\begin{equation*}
\begin{array}{|l||ccccccccccc|}\hline
|w|       & 1 & 2 & 3 & 4  & 5  & 6  & 7   & 8   & 9   & 10  & 11   \\ \hline
\mbox{maximal}\ |v| & 2 & 7 & 7 & 36 & 36 & 37 & 173 & 172 & 171 & 170 & 1230 \\ \hline
\end{array}
\end{equation*}
So, at least all words of length less than $12$ must occur with bounded gaps in the classical Kolakoski sequence supporting the conjecture that it is repetitive. Note that making this observation precise would prove that the Kolakoski sequence is repetitive, because this list tells us that there is no $C^{\infty}$-word of length greater than $2\times 11+1230=1252$ such that its prefix of length (less than) $11$ does not occur again within this word.
The jumps in this list are closely related to the ``degree'' that we introduce in the next section.

\section{\texorpdfstring{$C^{\infty}$-}{Smooth }words and the ``Kolakoski measure''}

We say that a $C^{\infty}$ has \emph{degree} $j$ if
\begin{equation*}
D^j(w)\neq\varepsilon, \qquad D^{j+1}(w)=\varepsilon.
\end{equation*}
We call $C^{\infty}$-words of degree $0$, i.e., the primitives of the empty word $\varepsilon$, \emph{fundamental words}. Note that a fundamental word has length less than $\max\{s,2r+1\}$.

We now define a function $\mu$ on the \emph{cylinder sets} $[w]$ of $\ABC^{\NN}$, i.e., $[w]=[w_1\ldots w_n]=\linebreak \{z\in\ABC^{\NN}\mathbin| z_1=w_1,\ldots,z_n=w_n\}$, by
\begin{equation*}
\mu([w])=\begin{cases} \mu\left(\left[D^j(w)\right]\right)\cdot \frac1{(r+s)^j} & \mbox{if $w$ is a $C^{\infty}$-word of degree $j$},\\
0 & \mbox{if $w$ is not a $C^{\infty}$-word}.
\end{cases}
\end{equation*}
Here, we have to fix the function $\mu$ for all fundamental words, and we do so by requiring that $\mu([w])=\mu([\tilde{w}])$, $\mu([w])=\mu([\stackrel{\leftarrow}{w}])$ for all fundamental words $w$ and that $\sum_{w\in\ABC^n} \mu([w]) = 1$ for $1\le n< \max\{s,2r+1\}$. For example, one has for the fundamental words\\[2mm]
\begin{tabular}{lllll}
using $\ABC=\{1,2\}$: & $\mu([1]) = \frac12$ & $\mu([2]) = \frac12$ 
& $\mu([12]) = \frac13$ & $\mu([21]) = \frac13$ \\[2mm] 
using $\ABC=\{2,3\}$: & $\mu([2]) = \frac12$ & $\mu([3]) = \frac12$ 
& $\mu([23]) = \frac15$ & $\mu([32]) = \frac15$ \\
& $\mu([22]) = \frac3{10}$ & $\mu([33]) = \frac3{10}$ 
& $\mu([223]) = \frac15$ & $\mu([332]) = \frac15$ \\
& $\mu([233]) = \frac15$ & $\mu([322]) = \frac15$ 
& $\mu([2233]) = \frac15$ & $\mu([3322]) = \frac15$
\end{tabular}

Clearly, one has the property $\mu([D(w)]) = (r+s)\cdot \mu([w])$ for any $C^{\infty}$-word of length greater than or equal to $\max\{s,2r+1\}$, and one can use this to show:

\begin{thm}
For any $\ABC=\{r,s\}$, the function $\mu$ extends to a Borel-measure (also denoted $\mu$) on $\ABC^{\NN}$. This measure is mirror invariant, reversal invariant and shift invariant.
\end{thm}

\noindent
{\it Proof.} 
A careful case study as in \cite[Theorem 5.1]{Dek97} also works in the general case. \hfill$\Box$

The aim of introducing this measure is to connect it somehow to the frequencies of subwords $w$ in a Kolakoski sequence. Indeed, one can show:

\begin{thm}
Suppose that $z$ is a Kolakoski sequence over $\ABC=\{r,s\}$, where one of the numbers $r,s\in\NN$ is odd and the other even, and that the frequencies $f_w=\lim\limits_{n\to\infty} \frac{|z_1\ldots z_n|_w}{n}$ exist and satisfy $f_w=f_{\tilde{w}}$ for all words occuring in $z$. Then for all words $w$ we have $f_w=\mu([w])$.
\end{thm}

\noindent
{\it Proof.} 
The proof of \cite[Proposition 5.1]{Dek97} carries over to the general case, see \cite[Proposition 2.5]{diplom}. \hfill$\Box$

This is nice result -- if only we would know that the frequencies satisfy the required properties. In fact, one can state Keane's question for all Kolakoski sequences where one of the letters is odd and other one is even:
\begin{quote}
In a Kolakoski sequence over $\{r,s\}$ (where one letter is odd and the other one is even), does the frequency of $r$ exist and if so, does it equal $\frac12$?
\end{quote}

Much computing time has been dedicated to find evidences for or against the conjecture that the letter freqeuncy is $\frac12$. The numerical evidences against it are usually dismissed by looking at larger and larger parts of the Kolakoski sequence, see~\cite{Ste06}.

Since already the existence of the letter frequency is in question, one can try to find bounds on $\limsup_{n\to\infty} |z_1\ldots z_n|_r/n$ and $\liminf_{n\to\infty} |z_1\ldots z_n|_r/n$ using the $C^{\infty}$-words. A brute force approach is, of course, to generate all $C^{\infty}$-words of a certain length, say $n$, and check for those with the least number\footnote
{
	I.e., we have $a=\min\{|w|_r\mathbin| |w|=n\mbox{ and w is a }C^{\infty}\mbox{-word}\}$ 
}
$a$ of $r$s (since for any $C^{\infty}$-word $w$, its mirrored version $\tilde{w}$ is also a $C^{\infty}$ word, the maximal number of $r$s is $n-a$).  One then has\footnote{
	For a proof see \cite[Section 3.2]{KRpre}.
}
\begin{equation*}
\frac{a}{n}\le \liminf_{n\to\infty}  \frac{|z_1\ldots z_n|_r}{n} \le \frac12 \le   \limsup_{n\to\infty}\frac{|z_1\ldots z_n|_r}{n} \le \frac{n-a}{n}.
\end{equation*}
For example, one finds the following numbers for alphabets with $r+s\le 7$:
\begin{equation*}
\begin{array}{|l||c|c|c|c|c|c|}\hline
\mbox{alphabet} & \{1,2\} & \{2,3\} & \{1,4\} & \{3,4\} & \{2,5\} & \{1,6\} \\ \hline\hline
\mbox{length }n & 1355 & 8003 & 1131 & 1000 & 1000 & 1000 \\
a=\min\limits_{|w|=n} |w|_r & 669 & 3989 & 511 & 493 & 481 & 451 \\
\mbox{letter freq.} &  {\scriptstyle  0.5\pm0.0063} & {\scriptstyle  0.5\pm0.0016} & {\scriptstyle  0.5\pm0.0482} & {\scriptstyle  0.5\pm0.007} & {\scriptstyle  0.5\pm0.019} & {\scriptstyle  0.5\pm0.049} \\ \hline 
\end{array}
\end{equation*}

Alternatively, one can use a generating function approach, see \cite{KRpre} (based on \cite{NZ99}): For each word $w$ on the alphabet $\{r,s\}$, one defines the its \emph{weight} as polynomial $x^{|w|_r}\,y^{|w|_s}\,t^{|w|}$. By summing these weights over all $C^{\infty}$-words\footnote{
	In fact, it is computationally more feasible to sum over all words that just avoid to be $C^{\infty}$-words, i.e., words on $\{r,s\}$ that are not $C^{\infty}$-words but any of its (genuine) subwords is. This is the method used in \cite{KRpre,NZ99}.
},
a lower bound on the frequency is obtained by looking at the minimal degree of $x$ for a given power $t^n$. The bound $\frac12\pm\frac{17}{762} \approx 0.5\pm 0.0223097$ was obtained using this method for the alphabet $\ABC=\{1,2\}$..

\section{Chvatal's bound on the letter frequency}\label{sec:Chvatal}

Instead of considering $C^{\infty}$-words, Chvatal \cite{Chv94} in his unpublished technical report looked at infinite words over $\{1,2\}$ with the property that their run-length sequence is also a sequence over the same alphabet. A sequence over $\{r,s\}$ is said to be $1$-special. If only runs of length $r$ and $s$ occur in this sequence, we say that the sequence is $2$-special. And if in this run-length sequence only runs of length $r$ and $s$ occur, we call the original sequence $3$-special. We continue in this way and note that a Kolakoski sequence is $d$-special for all $d\in\NN$. 

We now write (a $d$-special) sequence and its $d$ iterated run-length sequences in a special way in an array: The first row is the original sequence, the first row its run-length sequence and so on, but we align them appropriately in the columns. E.g., for the classical Kolakoski sequence (here we use the Kolakoski sequence starting with $1$), we write 
\begin{equation*}
\begin{array}{ccccccccccccccccccccc}
1 & 2 & 2 & 1 & 1 & 2 & 1 & 2 & 2 & 1 & 2 & 2 & 1 & 1 & 2 & 1 & 1 & 2 & 2 & 1 & \ldots \\ 
1 &    & 2 &    & 2 & 1 & 1 &    & 2 & 1 &    & 2 &    & 2 & 1 &    & 2 &    & 2 & 1 & \ldots \\  
1 &    &    &    & 2 &    & 2 &    & 1 & 1 &    &    &    & 2 & 1 &    &    &    & 2 &    & \ldots \\
1 &    &    &    &    &    & 2 &    &    & 2 &    &    &    & 1 & 1 &    &    &    &    &    & \ldots \\
\end{array}
\end{equation*}
We now call the $i$th element in a the original sequence \emph{$d$-special} if the sequence itself is $d$-special and the $i$th column in this array has length at least $d$ (there are no blanks in the first $d$ lines of this column). We call this column (of length $d$) the \emph{type} of the corresponding $d$-special element in the sequence. E.g., the third letter in the Kolakoski sequence above is $2$-special of type $22$, while the $7$th letter is $4$-special of type $1122$ (it is also $2$-special of type $11$).

Now, the observation is that the type of a $d$-special element determines the first $d$ terms of the type of the previous $d$-special element as well as all the letters between them, and the type of a $d$-special element and the last term of the type of the next $d$-special element determine the remaining terms of the type of this next $d$-special element. These properties can be used to iterative build graphs $G_d$, $d\ge 1$. Since the same observations can be made about Kolakoski sequence on any alphabet $\{r,s\}$, we describe the more general case here.

The vertices of the graph $G_d$ are the types of the $d$-special elements. Since all elements of $\ABC^d$ occur as types, the graph $G_d$ has $2^d$ vertices.
We connect two vertices $u$ and $v$ by a directed edge $u\stackrel{ w}{\longrightarrow}v$ labelled $ w$, if $v$ is the next $d$-special type after $u$ in a $d$-special sequence $z$. If $v$ is the $i$th element of such a $d$-special sequence and $u$ the $j$th element, then the label $ w$ is the word $z_{i+1}z_{i+2}\ldots z_{j-1}z_j$. So if we follow any (infinite) directed path in such a graph and read the edge-labels, we get a $d$-special sequence. Conversely, any $d$-special sequence arises as such an infinite path. 

The trick is that one can build the graph $G_{d+1}$ from the graph $G_d$: 
\begin{itemize}
\item A path $Ar\ \stackrel{ w_1}{\longrightarrow}\ B_1s\ \stackrel{ w_2}{\longrightarrow}\ B_2s\ \stackrel{ w_3}{\longrightarrow}\ \ldots \ \stackrel{ w_r}{\longrightarrow}\ B_rs$ in the graph $G_d$ gives rise to the edges $Arr\stackrel{ w_1 w_2\ldots w_r}{\longrightarrow}B_rsr$ and $Ars\stackrel{ w_1 w_2\ldots w_r}{\longrightarrow}B_rsr$ in $G_{d+1}$. 
\item A path $As\ \stackrel{ w_1}{\longrightarrow}\ B_1r\ \stackrel{ w_2}{\longrightarrow}\ B_2r\ \stackrel{ w_3}{\longrightarrow}\ \ldots \ \stackrel{ w_r}{\longrightarrow}\ B_rr$ in the graph $G_d$ gives rise to the edges $Asr\stackrel{ w_1 w_2\ldots w_r}{\longrightarrow}B_rrr$ and $Ass\stackrel{ w_1 w_2\ldots w_r}{\longrightarrow}B_rrr$ in $G_{d+1}$.
\item A path $Ar\ \stackrel{ w_1}{\longrightarrow}\ B_1s\ \stackrel{ w_2}{\longrightarrow}\ B_2s\ \stackrel{ w_3}{\longrightarrow}\ \ldots \ \stackrel{ w_s}{\longrightarrow}\ B_ss$ in the graph $G_d$ gives rise to the edges $Arr\stackrel{ w_1 w_2\ldots w_r}{\longrightarrow}B_sss$ and $Ars\stackrel{ w_1 w_2\ldots w_r}{\longrightarrow}B_sss$ in $G_{d+1}$.
\item A path $As\ \stackrel{ w_1}{\longrightarrow}\ B_1r\ \stackrel{ w_2}{\longrightarrow}\ B_2r\ \stackrel{ w_3}{\longrightarrow}\ \ldots \ \stackrel{ w_s}{\longrightarrow}\ B_sr$ in the graph $G_d$ gives rise to the edges $Asr\stackrel{ w_1 w_2\ldots w_s}{\longrightarrow}B_srs$ and $Ass\stackrel{ w_1 w_2\ldots w_s}{\longrightarrow}B_srs$ in $G_{d+1}$.
\end{itemize}

The graphs $G_1$ and $G_2$ are :
\begin{equation*}
\xymatrix{
\\
*+[F]{r} \ar@(ul,dl)[]|{r} \ar@/^/[rr]|s
&& *+[F]{s} \ar@(dr,ur)[]|{s} \ar@/^/[ll]|{r} \\
}\qquad\qquad
\xymatrix{
*+[F]{rr} \ar@/^/[dd]|{s^s} \ar@/^/[rr]|{s^r}
&& *+[F]{sr} \ar@/^/[dd]|{r^s} \ar@/^/[ll]|{r^r} \\ \\
*+[F]{ss} \ar@/^/[uu]|{r^r} \ar@/^/[rr]|{r^s}
&& *+[F]{rs} \ar@/^/[uu]|{s^r} \ar@/^/[ll]|{s^s}
}
\end{equation*}
To get bounds on the letter frequencies from a graph $G_d$, one associates to an edge with edge label $w$ the cost $x\cdot |w|-|w|_r$. If one now uses for $x$ a number $\frac12\le x<1$ that is smaller than the maximal possible letter frequency that can occur for a $d$-special sequence, then one finds a negative cycle in this graph. Applying this method to $G_6$ for alphabets with $r+s\le 7$, one finds the following bounds:
\begin{equation*}
\begin{array}{|l||c|c|c|c|c|c|}\hline
\mbox{alphabet} & \{1,2\} & \{2,3\} & \{1,4\} & \{3,4\} & \{2,5\} & \{1,6\} \\ \hline\hline
\mbox{upper bound} & {}^{12}\!/\!_{23} & {}^{53}\!/\!_{105} & {}^{592}\!/\!_{1085} & {}^{46}\!/\!_{91} & {}^{4834}\!/\!_{9527} & {}^{1478}\!/\!_{2821} \\
\mbox{letter freq.} &  {\scriptstyle  0.5\pm0.0218} & {\scriptstyle  0.5\pm0.0048} & {\scriptstyle  0.5\pm0.0457} & {\scriptstyle  0.5\pm0.0055} & {\scriptstyle  0.5\pm0.0075} & {\scriptstyle  0.5\pm0.0240} \\ \hline 
\end{array}
\end{equation*}
By a clever use of the structure of the graphs $G_d$ and efficient programming, Chvatal used $G_{22}$ in \cite{Chv94} which yields the upper bound $616904/1231743$ for the classical Kolakoski sequence over $\ABC=\{1,2\}$, i.e., the letter frequencies are confined to $0.5\pm0.000838$.

\section{Squares (and Cubes)}

The question which (and how many) squares occur in the classical Kolakoski sequence was asked by \cite{Pau93}. Shortly thereafter, Carpi \cite{Car93,Car94} and Lepist\"{o} \cite{Lep94} answered the question by finding all squares the occur: in the classical Kolakoski sequence (i.e., using the alphabet $\ABC=\{1,2\}$) only squares of length $1$, $2$, $3$, $9$ and $27$ (\cite[Theorem 1]{Lep94}, \cite[Proposition 1]{Car93}, \cite[Proposition 3]{Car94}) occur; in particular, it is cube-free (\cite[Corollary 1]{Lep94},\cite[Proposition 2]{Car93}, \cite[Proposition 4]{Car94}). Here, a square $w$ of length $n$ is a $C^{\infty}$-word with $|w|=n$ such that $ww$ is also a $C^{\infty}$-word.

\begin{table}[t]
{
\begin{displaymath}
\begin{array}{|c||c|c|c|}\hline
\rule[-2mm]{0cm}{.5cm}
\mbox{length }n & {\begin{array}{c} \mbox{number of}\\ \mbox{squares}\end{array}} &
   {\begin{array}{c} \mbox{complexity}\\ \gamma(n)\end{array}} & 
   {\begin{array}{c} \mbox{max.} \\ \mbox{MOOR}\end{array}} \\ \hline \hline  
\rule[-1.5mm]{0cm}{.35cm}
1 & 2 & 2 & 2 \\
\rule[-1.5mm]{0cm}{.3cm}
2 & 2 & 4 & 2 \\
\rule[-1.5mm]{0cm}{.3cm}
3 & 6 & 6 & 2 \,{}^2\!/\!_3 \\
\rule[-1.5mm]{0cm}{.3cm}
9 & 12 & 42 & 2 \,{}^1\!/\!_9 \\
\rule[-1.5mm]{0cm}{.3cm}
27 & 24 & 486 & 2 \,{}^1\!/\!_{27} \\ \hline 
\end{array}
\end{displaymath}}
\caption{The classical Kolakoski case over $\ABC=\{1,2\}$: There are $46$ squares and no cubes.\label{tab:12}}
\end{table}

\begin{table}[t]
{
\begin{displaymath}
\begin{array}{|c||cc|c|c|}\hline
\rule[-2mm]{0cm}{.5cm}
\mbox{length }n & {\begin{array}{c} \mbox{number of}\\ \mbox{squares}\end{array}} &
   {\begin{array}{c} \mbox{thereof also}\\ \mbox{cubes}\end{array}} & 
   {\begin{array}{c} \mbox{complexity}\\ \gamma(n)\end{array}} & 
   {\begin{array}{c} \mbox{max.} \\ \mbox{MOOR}\end{array}} \\ \hline \hline  
\rule[-1.5mm]{0cm}{.35cm}
1 & 2 & 2 & 2 &  3 \\
\rule[-1.5mm]{0cm}{.3cm}
4 & 4 &   & 8 & 2 \,{}^1\!/\!_2 \\
\rule[-1.5mm]{0cm}{.3cm}
6 & 4 &   & 14 & 2 \,{}^1\!/\!_6 \\
\rule[-1.5mm]{0cm}{.3cm}
10 & 20 & 2 & 30 & 3 \\
\rule[-1.5mm]{0cm}{.3cm}
15 & 30 &   & 58 &  2 \,{}^7\!/\!_{15} \\
\rule[-1.5mm]{0cm}{.3cm}
25 & 100 & 16 & 130 & 3 \,{}^4\!/\!_{25} \\ \hline
\end{array}
\end{displaymath}}
\caption{For $\ABC=\{2,3\}$, there are $160$ squares and $20$
  cubes among the $C^{\infty}$-words.\label{tab:23}}
\end{table}

\begin{table}[!th]
{
\begin{displaymath}
\begin{array}{|c||ccc|c|c|}\hline
\rule[-2mm]{0cm}{.5cm}
\mbox{length }n & {\begin{array}{c} \mbox{number of}\\ \mbox{squares}\end{array}} &
   {\begin{array}{c} \mbox{thereof also}\\ \mbox{cubes}\end{array}} &  {\begin{array}{c} \mbox{thereof also}\\
  \mbox{fourth powers}\end{array}} & {\begin{array}{c} \mbox{complexity}\\ 
  \gamma(n)\end{array}} & 
   {\begin{array}{c} \mbox{max.} \\ \mbox{MOOR}\end{array}} \\ \hline \hline  
\rule[-1.5mm]{0cm}{.35cm}
1 & 2 & 2 & 2 & 2 &  4 \\
\rule[-1.5mm]{0cm}{.35cm}
2 & 4 & 2 &   & 4 &  3 \\
\rule[-1.5mm]{0cm}{.35cm}
5 & 10 & 8 &   & 20 &  3 \,{}^2\!/\!_5 \\
\rule[-1.5mm]{0cm}{.35cm}
8 & 6 & & & 36 & 2 \,{}^1\!/\!_4 \\
\rule[-1.5mm]{0cm}{.35cm}
13 & 4 & & & 96 &  2 \,{}^1\!/\!_{13} \\
\rule[-1.5mm]{0cm}{.35cm}
20 & 30 & & & 198 &  2 \,{}^7\!/\!_{20} \\
\rule[-1.5mm]{0cm}{.35cm}
40 & 8 & & & 630 & 2 \,{}^1\!/\!_{40} \\
\rule[-1.5mm]{0cm}{.35cm}
50 & 152 & & & 964 &  2 \,{}^{12}\!/\!_{25} \\
\rule[-1.5mm]{0cm}{.35cm}
100 & 48 & & & 3124 &  2 \,{}^3\!/\!_{50} \\
\rule[-1.5mm]{0cm}{.35cm}
116 & 364 & & & 4160 &  2 \,{}^{59}\!/\!_{116} \\
\rule[-1.5mm]{0cm}{.35cm}
134 & 400 & & & 5438 &  2 \,{}^{65}\!/\!_{134} \\
\rule[-1.5mm]{0cm}{.35cm}
174 & 8 & & & 8658 & 2 \,{}^1\!/\!_{174} \\
\rule[-1.5mm]{0cm}{.35cm}
241 & 144 & & & 14694 &  2 \,{}^{20}\!/\!_{241} \\
\rule[-1.5mm]{0cm}{.35cm}
259 & 100 & & & 16588 &  2 \,{}^2\!/\!_{37} \\
\rule[-1.5mm]{0cm}{.35cm}
272 & 864 & & & 18358 &  2 \,{}^{145}\!/\!_{272} \\
\rule[-1.5mm]{0cm}{.35cm}
308 & 960 & & & 23288 &  2 \,{}^1\!/\!_2 \\
\rule[-1.5mm]{0cm}{.35cm}
317 & 960 & & & 24554 &  2 \,{}^{160}\!/\!_{317} \\
\rule[-1.5mm]{0cm}{.35cm}
353 & 1044 & & & 29738 & 2 \,{}^{169}\!/\!_{353} \\
\rule[-1.5mm]{0cm}{.35cm}
408 & 4 & & & 38462 & 2 \\
\rule[-1.5mm]{0cm}{.35cm}
417 & 16 & & & 40046 & 2 \,{}^1\!/\!_{139} \\
\rule[-1.5mm]{0cm}{.35cm}
453 & 28 & & & 46474 &  2 \,{}^2\!/\!_{151} \\
\rule[-1.5mm]{0cm}{.35cm}
644 & 2072 & & & 82292 &  2 \,{}^{177}\!/\!_{322} \\
\rule[-1.5mm]{0cm}{.35cm}
716 & 2252 & & & 100990 &  2 \,{}^{375}\!/\!_{716} \\
\rule[-1.5mm]{0cm}{.35cm}
734 & 2312 & & & 106410 &  2 \,{}^{375}\!/\!_{734} \\
\rule[-1.5mm]{0cm}{.35cm}
806 & 2492 & & & 126570 &  2 \,{}^{399}\!/\!_{806} \\
\rule[-1.5mm]{0cm}{.35cm}
975 & 12 & & & 177330 &  2 \,{}^2\!/\!_{975} \\
\rule[-1.5mm]{0cm}{.35cm}
1065 & 12 & & & 208018 &  2 \,{}^2\!/\!_{1065} \\
\rule[-1.5mm]{0cm}{.35cm}
1529 & 4960 & & & 376874 &  2 \,{}^{845}\!/\!_{1529} \\
\rule[-1.5mm]{0cm}{.35cm}
1691 & 5404 & & & 451208 &  2 \,{}^{929}\!/\!_{1691} \\
\rule[-1.5mm]{0cm}{.35cm}
1709 & 5404 & & & 460688 &  2 \,{}^{896}\!/\!_{1709} \\
\rule[-1.5mm]{0cm}{.35cm}
1745 & 5500 & & & 480304 &  2 \,{}^{178}\!/\!_{349} \\
\rule[-1.5mm]{0cm}{.35cm}
1871 & 5860 & & & 550730 &  2 \,{}^{983}\!/\!_{1871} \\
\rule[-1.5mm]{0cm}{.35cm}
1925 & 12020 & & & 581470 &  2 \,{}^{989}\!/\!_{1925} \\
\rule[-1.5mm]{0cm}{.35cm}
2105 & 6508 & & & 684994 &  2 \,{}^{1049}\!/\!_{2105} \\ \hline
\end{array}
\end{displaymath}}
\caption{For $\ABC=\{1,4\}$, there are $59\,964$ squares of which are $12$ also cubes and 
  only $2$ are also fourth powers among the $C^{\infty}$-words.\label{tab:14}}
\end{table} 

The algorithm for finding squares is based on the following observations: If $ww$ is a square, its derivative has the form $D(ww) = uvu$ where $|uv|$ has to be even (otherwise, not $ww$ but $w\tilde{w}$ will be a primitive) and $|v|\le 1$ (we have $D(w)=u$ and $v$ arises because of the rule on how to derive first and last runs). There is one speciality, though, if $r<\frac{s}2$: In these cases, $v$ might also be a ``negative'' power $r^{-1}$ or $s^{-1}$ of length $-1$, meaning that in $uv$ the $v$ cancels the last symbol of $u$. 

Now if one continues differentiating, one gets a sequence of words $D(ww)=u_1v_1u_1$, $D^2(ww)=u_2v_2u_2$, \ldots $D^k(ww)=u_kv_ku_k$. But one can show that in this sequence the length of $|v_i|$ is bounded by $-1\le|v_k|\le 2s+1$ where the lower bound $-1$ only can appear if $r<\frac{s}2$, see \cite[Lemma 4.4]{diplom} (compare to \cite[Lemma 1]{Lep94} for $\ABC=\{1,2\}$). Furthermore, we must always have that $|u_iv_i|$ is even for all $i$. Thus, one now has an algorithm to find squares in Kolakoski sequences: We start with all $C^{\infty}$-words of the form $uvu$ where $u$ is a fundamental word, $-1\le|v|\le 2s+1$ and $|uv|$ is even (and/or $v=\varepsilon$ and thus we already have a square $uu$). Then construct all primitives which are again of the form $u'v'u'$ with either $|u'v'|$ even and $-1\le|v'|\le 2s+1$, and/or which happen to be a square $u'u'$. Continue in this way. If there are eventually no more words of this form left, the algorithm stops and one has calculated all squares among the $C^{\infty}$-words.

We note, however, that it is a priori not clear whether this algorithm will indeed stop, or if there are only finitely many squares in a Kolakoski sequence besides the classical one. However, we used this algortihm to check the $C^{\infty}$-words over the alphabets $\{2,3\}$ and $\{1,4\}$: Both, similar to the classical Kolakoski sequence, have only finitely many squares -- there are a total of $160$ different squares of smooth words over $\{2,3\}$ but $59\,964$ squares of smooth words in $\{1,4\}$. We list the numbers in these case together with the classical Kolakoski sequence in Tables~\ref{tab:12}--\ref{tab:14}. We also listed the number of cubes and fourth powers in this cases together with the complexity at this word length. The \emph{maximal order of repetition}, for short \emph{MOOR}, or \emph{repetition exponent} for a $C^{\infty}$-word $w=uv$ is given by the maximum of $|ww\ldots wu|/|w|$ such that $ww\ldots wu$ is also a $C^{\infty}$-word. 

Since there are only finitely many squares, the corresponding Kolakoski sequences cannot be obtained by a (usual)  
 substitution rule, see \cite[Theorem 2]{Lep94} (if a substitution sequence has one square, one gets infinitely many using the substitution rule repeatedly). Also, the following conjecture was stated in \cite{Car94}: For any repetition exponent $q>1$, the length of $C^{\infty}$-words having this exponent is bounded.

\section{Palindromes}

If a $C^{\infty}$-word is a palindrome, i.e., if we have $w=\stackrel{\leftarrow}{w}$, then $D(w)$ is also a palindrome. Conversely, only palindromes of odd length have primitives that are also palindromes. We have a look at the following table: 
\begin{displaymath}
\begin{array}{|c|c|}\hline
\mbox{palindrome} & \mbox{primitives} \\ \hline \hline
22  & 1122,\, 21122,\, 11221,\, 211221,\\
    & 2211,\, 12211,\, 22112,\, 122112 \\ \hline
212 & \boldsymbol{11211},\, 211211,\, 112112,\, \boldsymbol{2112112}, \\
    & \boldsymbol{22122},\, 122122,\, 221221,\, \boldsymbol{1221221} \\ \hline
121 & \boldsymbol{121121},\, \boldsymbol{212212} \\ \hline
\end{array}
\end{displaymath}
Thus, together with a further observation, one has in fact an algorithm how to construct palindromes, compare \cite[Section 4.1.3]{Lad99} (see also \cite{BL03,BDLV05}): Start with all palindromic fundamental words. Palindromes of odd length where the letter in the middle is odd, will have palindromes of odd length among their primitives. Palindromes of odd length where the letter in the middle is even, will have palindromes of even length among their primitives. Palindromes of even length do not have palindromic primitives.

Since $w$ is a palindrome iff $\tilde{w}$ is a palindrome, palindromes (in fact, palindromic fundamental words) of odd length where the letter in the middle is odd play a special role and can be used to construct all palindromes. For example by repeatedly constructing primitives, one gets the following palindromic two-sided infinite sequence with the number $1$ ``in the middle'' when starting if the fundamental word $1$ over $\ABC=\{1,2\}$:
\begin{displaymath}
\begin{array}{r|c|l}
\ldots122121122 & 1 & 221121221\ldots \\
\end{array}
\end{displaymath}
Applying the operation $\tilde{\cdot}$ to this word yields:    
\begin{displaymath}
\begin{array}{r|c|l}
\ldots211212211 & 2 & 112212112\ldots \\
\end{array}
\end{displaymath}
The primitives of this infinite sequence are:
 \begin{displaymath}
\begin{array}{r|c|l}
\ldots12112212 & 11 & 21221121\ldots \\
\ldots21221121 & 22 & 12112212\ldots \\
\end{array}
\end{displaymath}
Consequently, looking at the symmetric part of these sequences, one has $2$ palindromes of each length (for details see the cited literature). In fact, one always has the single letters as fundamental words, there are for all alphabets at least $2$ palindromic $C^{\infty}$-words for each length. E.g., for $\ABC=\{2,3\}$, the same construction as before works, where we now have
\begin{displaymath}
\ldots 2223322332223|3|3222332233222 \ldots
\end{displaymath}
The situation gets a bit more complicated if there is more than one palindromic fundamental word of odd length with an odd letter in the middle. E.g., the alphabet $\ABC=\{1,4\}$ the two fundamental words with the stated property are $1$ and\footnote
{
	Note that $111$ is a single run of length $3 < 4=s$ and thus we have $D(111)=\varepsilon$.
} 
$111$. Thus, additional palindromes appear (but for each length, one has at most two times the number of palindromic fundamental words of odd length with odd letter in its middle):
\begin{displaymath}
\begin{array}{|c|l|}\hline
\mbox{length} & \mbox{palindromes} \\ \hline \hline
1 & 1,\, 4 \\
2 & 11,\, 44 \\
3 & 111,\, 414,\, 444,\, 141 \\
4 & 1111,\, 4444 \\
5 & 44144,\, 14141,\, 11411,\, 41414 \\
6 & 411114,\, 144441 \\
7 & 4441444,\, 1114111 \\
8 & 14111141,\, 44111144,\, 41444414,\, 11444411 \\ \hline
\end{array}
\end{displaymath}

Generalizations of palindromes, namely words of the form $\stackrel{\leftarrow}{w}\!\!vw$ (``palindromes with a gap in the middle'') have been studied in \cite{Hua08,Hua09}.

\section{Complexity}

It is clear that the set of subwords of a Kolakoski sequence is a subset of the $C^{\infty}$-words over the same alphabet. Since one, in fact, conjectures that the two sets are even identical, one tries to establish bounds on the number of $C^{\infty}$-words for a given length. We denote the complexity of $C^{\infty}$-words, i.e., the number of $C^{\infty}$-words of length $n$, by $\gamma(n)$. 

Again, one can straightforwardly generalize results by Dekking.

\begin{thm}
Let $\gamma(n)$ be the number of $C^{\infty}$-words of length $n$ in the alphabet $\ABC=\{r,s\}$. Then
\begin{enumerate}
\item there is an $N\in\NN$ such that $\gamma(n) \le n^{\alpha}$ where $\alpha = \frac{\ln (2
s^2)}{\ln (\frac{2rs}{r+s})}$ for all $n \ge N$.
\item there is an $N\in\NN$ and a constant $C>0$ such that $\gamma(n) \ge C \cdot n^{\beta}$ where $\beta = \frac{\ln(r+s)}{\ln (\frac{r^2+s^2}{r+s})}$ for all $n \ge N$.
\end{enumerate}
\end{thm}

\noindent
{\it Proof.}
For a proof in the classical case $\ABC=\{1,2\}$, see \cite[Propositions 3 \& 4]{Dek80}. For the generalizations, see
 \cite[Propositions 4.1 \& 4.3]{diplom}.\hfill$\Box$
 
For the alphabet $\ABC=\{1,2\}$ these bounds have recently been improved in \cite{HW10} (based on previous work \cite{Wea89}). In this case, there are positive constants $C_1,C_2$ such that 
\begin{equation*}
C_1\,n^{2.7087}< \gamma(n)< C_2\,n^{2.7102}
\end{equation*}
for all $n\in\NN$.

\begin{figure}[t]
\centerline{\includegraphics[width=.6\textwidth]{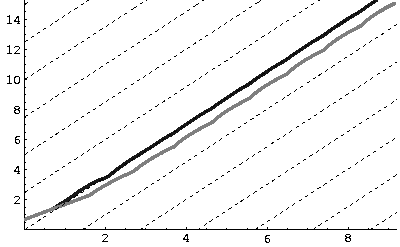}}
\caption{$\ln\gamma(n)$ vs.\ $\ln n$ for $\ABC=\{1,4\}$ (dark gray) and $\ABC=\{2,3\}$ (light gray). The dotted lines are the graphs of $f_k(n)=n^{k\cdot\ln 5/\ln\frac52}$ in this double-log plot.\label{fig:complexity}}
\end{figure}

In fact, one can conjecture:
\begin{quote}
There are positive constants $C_1, C_2$ such that 
\begin{displaymath}
C_1 \cdot n^{\delta} \le \gamma(n) \le C_2 \cdot n^{\delta}, \mbox{ where } \delta =\frac{\ln
  (r+s)}{\ln \frac{r+s}{2}}. 
\end{displaymath} 
\end{quote}
Noting that for $\ABC=\{1,2\}$ we have $\delta=\ln3/\ln\frac32\approx 2.7095$, we see that this conjecture is well supported by the above result. For $\ABC=\{2,3\}$ and $\ABC=\{1,4\}$, we refer to numerical results that we show in Fig.~\ref{fig:complexity}.

\bibliographystyle{amsplain}

\end{document}